\documentclass{amsart}[12pt]
\usepackage{amscd}

\usepackage{amsmath}
\usepackage{amssymb}
\newcommand{\R}{{\mathbb R}}
\newcommand{\C}{{\mathbb C}}
\newcommand{\Z}{{\mathbb Z}}

\newcommand{\ah}{\mathcal{A}_H}
\newcommand{\rs}{{\mathbb R} \times S^1}

\newcommand{\hl}{H_{\lambda}}
\newtheorem{theo}{Theorem}
\newtheorem{lm}{Lemma}
\newtheorem{df}{Definition}
\newtheorem{prop}{Proposition}
\newtheorem{cor}{Corollary}
\newtheorem{rema}{Remark}

\begin{document}

\title[A note on the representability]{A note on the representability of a certain Hamiltonian capacity}

\author{Dragomir L. Dragnev\\}

%\address{Courant Institute of Mathematical Sciences\\
%New York University \\ 251 Mercer street \\New York, NY 10012}
%\email{dd22@nyu.edu}
%\thanks{Preliminary version.}
\begin{abstract}
In this note we establish a representation property for a certain
Hamiltonian capacity on $\R^{2n}$ with the standard symplectic
structure. We demonstrate that the value of this capacity on an
open set with a contact type boundary is an element of the action
spectrum of the boundary.
\end{abstract}

\maketitle
%%%%%%%%%%%%%%%%%%%%%%%%%%%%%%%%%%%%%%%%%%%%%%%%%%%%%%%%%%%%%%%%%%%%%%%%%%%
\section{Introduction}
Consider $\R^{2n}$ with the standard symplectic structure
$\omega_0 = \sum_{j=1}^n dx_j \wedge dy_j$. Recall the definition
of a symplectic capacity on $(\R^{2n} \simeq \C^n, \omega_0)$.
\begin{df}\label{capacity} A (normalized) symplectic capacity is a map which associates to a
given set $U \subset \C^n$ a number $c(U)$  with the following
properties,
\begin{enumerate}
\item \emph{Monotonicity:} If $U \subset V$ then $c(U) \leq c(V)$,
\item \emph{Symplectic invariance:} $c(\phi(U)) = c(U)$, for any
sympectomorphism $\phi$ of $\C^n$, \item \emph{Homogeneity:} $c(
aU) = a^2 c(U)$ for any real number $a$. \item
\emph{Normalization:} $c (B^{2n} (1)) = c (Z(1)) = \pi$, where
$B^{2n} (1)$ is the unit ball in $\C^n$, centered at the origin
and $Z(1) = \{ z=(z_1, \ldots , z_n) \in \C^n \mid |z_1| < 1 \}$
\end{enumerate}
\end{df}
Notice that it is sufficient to define such a map $c$ with the above
properties on open and bounded subsets of $\C^n$, afterwards we
can extend it to any open set as follows,
$$c(U) = \sup \{ c(V) \mid V  \textmd{ is bounded and connected and }
V \subset U \}$$ and to any subset by:
$$c(E) = \inf \{ c(U) \mid U \textmd{ is open and } E \subset U
\}$$ The notion ``symplectic capacity" was introduced by Ekeland
and Hofer, \cite{EH2,EH3}, although the first symplectic
invariant, satisfying the axioms of Definition \ref{capacity}, was
presented by Gromov in \cite{MG}. Now there are many capacity
functions derived from various constructions in symplectic
topology and for these we refer to the survey paper, \cite{CHLS}.

Consider a bounded domain $U \subset \R^{2n} \simeq \C^n$ with
smooth boundary $S= \partial U$. The restriction of the symplectic
form $\omega_0$ on $S$ gives rise to a line bundle over $S$,
$\mathcal{L}_{S} = \ker \omega_0|_{S} \subset TS$ that is,
$$\mathcal{L}_{S} =\{ (x, \xi) | \omega_0 (\xi,\eta) = 0 \textmd{ for all } \eta \in T_x
S\}$$ The integral curves of the line bundle $\mathcal{L}_{S}$,
are called \emph{characteristics} of $S$. Let $\alpha$ be any
primitive of $\omega_0$ on $\C^n$. If $\gamma$ is a closed
characteristic on $S$ we define its \emph{action} $A(\gamma)=
\int_{\gamma} \alpha$ and we define the \emph{action spectrum} of
$S$ to be the set:
$$\Sigma(S) = \{ k |A(\gamma)| | k \in \mathbb{N}, \gamma \textmd{ is a closed characteristic of }S
\}$$ We point out that if $S$ is a regular level surface of
autonomous Hamiltonian function $H =  \textmd{const}$, then the
Hamiltonian vector field $X_H$ is a section of $\mathcal{L}_{S}$
and so in this case the closed characteristics of $S$ are the
periodic orbits of $X_H$ on $S$. The problem of the existence of
closed characteristics on a hypersurface of $\R^{2n}$ has been
studied extensively and it depends on the symplectic properties of
$S$. Indeed results of M. Hermann, \cite{MH} and V. Ginzburg,
\cite{G1}, show that one should not expect the existence of closed
characteristics on an arbitrary hypersurface. On the other hand,
the validity of Weinstein's conjecture in $\R^{2n}$, proved by C.
Viterbo, \cite{V0}, guarantees existence of closed characteristics
provided that the hypersurface $S$ is of \emph{contact type}. We
recall that a hypersurface $S$ in $\C^n$ is called \emph{contact},
if there is a Liouville vector field $X$, i.e. $L_X \omega_0 =
\omega_0$, which is defined in a neighborhood of $S$ and is
transversal to $S$. The hypersurface $S$ is called of
\emph{restricted contact type} if the vector field $X$ is globally
defined on $\C^n$. We remark that S. Bates, \cite{Ba} found
examples of hypersurfaces which are of contact type, but not of
restricted contact type.

There are several interesting questions to be answered regarding
the capacities and one of them concerns their
\emph{representability} that is when $c(U) \in \Sigma(\partial
U)$? There are several known results in this direction. Most of
them concern open sets with restricted contact type (RCT)
boundary. It is known that Ekeland-Hofer capacities,
\cite{EH2,EH3}, Hofer-Zehnder capacity, \cite{HZ}, Viterbo
capacity, \cite{V2}, the Floer-Hofer capacity and etc. enjoy this
property for RCT open sets. In general the question of representability of the capacities is a delicate one. D. Hermann showed in \cite{He2}, that
not all symplectic capacities are representable even if $U$ is a
RCT open set. In a recent paper, \cite{ABHS}, the main result implies the existence of dynamically convex star-shaped levels in $\R^4$ where the Gromov width is not representable as the action of a closed characteristic. On the other hand all of the above mentioned
representable capacities fall in the class of Hamiltonian
capacities in D. Hermann's termionology, \cite{He2}. In the same
reference, \cite{He2}, Question 1.2.2, he poses the question if
all Hamiltonian capacities are representable on open sets with
contact type boundary?

In this note we give positive answer to this question for the
Floer-Hofer capacity in the terminology \cite{He, He1,D}. It falls
in the class of Hamiltonian capacities. Our main theorem is:
\begin{theo}\label{D} Let $c$ be the Floer-Hofer capacity and $U$ be a bounded domain in $\R^{2n}$ with
contact type boundary. Then $c(U) \in \Sigma(\partial U)$.
\end{theo}
We mention that there are similar representation theorems
concerning the Viterbo capacity, and Ekeland-Hofer capacity, for
particular open sets with contact type boundary, e.g. $U$ is a
tubular neighborhood of Lagrangian torus or a tubular neighborhood
of hyperbolic Lagrangian submanifold, cf. \cite{V3,Ba,Th}.

As a consequence of the above theorem, we can strengthen a little bit Theorem 1.3 in view of Remark 1.4 in \cite{D}. We have,
\begin{cor}\label{D22}
Let $S$ be a compact hypersurface in $(\R^{2n},
\omega_0)$ of contact type. Let $\phi$ be the time-1 map of a
compactly supported Hamiltonian $H$ on $[0,1] \times \R^{2n}$ such
that $E(\phi) \leq c (S)$. Then there exists $x \in S$ such that
$\phi(x) \in \mathcal{L}_x S$.
\end{cor}

Next, we provide a brief description of the
strategy of the proof of Theorem \ref{D}. It exploits ideas
similar to the ones used in \cite{V0,HZ1,EH2}. We thicken the
boundary $S=\partial U$, i.e. using the contact definition we
foliate a neighborhood of $S$, into diffeomorphic images of $S$,
$\{S_{\delta}\}$, $\delta \in (0, \epsilon)$ for small $\epsilon >
0$ and $S=S_0$. Then we find a sequence of Hamiltonians
$H_{\delta_k} \in \mathcal{H}_{ad} (U)$, (defined in the next section), and a sequence of
critical points $\{x_{\delta_k}\}$ of
$\mathcal{A}_{H_{\delta_k}}$, with the property that each
$x_{\delta_k}$ is a closed characteristic of $S_{\delta_k}$ and
\begin{equation}\label{cap1}
\sigma(H_{\delta_k}) = \mathcal{A}_{H_{\delta_k}} (x_{\delta_k})
\end{equation}
where $\sigma$ is an action selector coming from certain type of Morse homology for the Hamiltonian action functional $\mathcal{A}_{H}$, (see Section \ref{FHc}).
Then under certain assumptions we would have that
\begin{equation}\label{cap2}
c(U) = \lim_{k \to \infty} \sigma(H_{\delta_k})
\end{equation}
Now, if the lengths of the closed characteristics $x_{\delta_k}$
are uniformly bounded, we could find a subsequence which of
$\{x_{\delta_k}\}$ converging to $x_0$ - a closed characteristic
of $S=S_0$. Assuming that $\lim_{k \to \infty} H_{\delta_k}
(x_{\delta_k}) = 0$, we would get using, (\ref{cap1},\ref{cap2}),
that
\begin{equation}\label{repr}
c(U) = A(x_0) \in \Sigma(S)
\end{equation}
This way, the crucial point is to establish the uniform
boundedness of the lengths of the closed characteristics
$\{x_{\delta_k}\}$. In case $U$ is RCT open this is a consequence
of the global definition of the Liouville vector field or
equivalently, cf. \cite{HZ}, the existence of a 1-form $\alpha_0$
on $\C^n$ so that $d\alpha_0 = \omega_0$ on $\C^n$ and $\alpha_0
\wedge (d\alpha_0) ^{n-1}$ is a volume form on $S$. In case $U$ is
with contact type boundary we establish the boundedness by
utilizing an elegant idea of Ph. Bolle, \cite{B}, which was
successfully applied by the author in \cite{D}. Namely using the fact
that each $x_{\delta_k}$ is a deformation of a constant solution
for $H_{\delta_k}$, an application of Stokes' theorem implies the
desired boundedness.
\subsection*{Acknowledgments}I would like to thank O. Cornea, D. Hermann, K. Honda, U. Hryniewicz and E.
Kerman for the interest in this paper and their critical remarks and suggestions.

\section{The Floer-Hofer capacity.} \label{FHc}
The Hamiltonian capacity we are concerned here may be obtained by
the so-called \emph{action selector method}, cf. \cite{HZ,FGS}. We
describe it briefly. The idea is to consider the set of
\emph{admissible} Hamiltonians for a given bounded domain $U$. In
general these are functions which are bounded from below, (or
above depending on the convention), on $U$, (usually by $0$),
which are allowed to grow rapidly near $\partial U$ and have
certain growth rate at infinity. We denote this set by
$\mathcal{H}_{ad} (U)$. For $H \in \mathcal{H}_{ad} (U)$, we
consider the Hamiltonian action functional $\ah (x)$, $x \in
C^{\infty} (S^1, \R^{2n})$,
\begin{equation} \label{haf}
\ah (x) = \int_{S^1} x^{*} \lambda_0 - \int_0 ^1 H(t, x(t)) dt
\end{equation}
where $\lambda_0$ is some primitive of $\omega_0$. Then one
``selects", $\sigma(H)$ - a critical value of $\ah(x_{0})$ where
$x_0$ is some ``topologically visible" periodic orbit of the
Hamiltonian vector field $X_H$ and the capacity of $U$ is defined
by,
\begin{equation}\label{capa}
c(U) = \inf \{ \sigma(H) | H \in \mathcal{H}_{ad} (U)\}
\end{equation}
We would like to outline how the action selector method we just
described, may be used to construct the Floer-Hofer capacity. For
the sake of brevity we will just present the important steps in
the construction, to fix the notation. We omit some of the details
for which, however, we present references. The Floer-Hofer
capacity is based on the computations of the symplectic homology,
\cite{FH1} for balls in \cite{FHW}.

Define the set of admissible Hamiltonians, $\mathcal{H}_{ad}(U)$
for an open bounded set $U \subset \C^n$,
\begin{alignat*}{4}
\mathcal{H}_{ad}(U) = \{H: S^1 \times \C^n \to \R| H(t,z) \leq 0
\textmd{ on } \bar{U} \textmd{ and } \\ H(t,z) = \mu |z|^2
\textmd{ outside of a compact set}\}
\end{alignat*}
The Floer homology, \cite{F,Sa}, can be thought as an
infinite-dimensional version of Morse theory for the Hamltonian
action functional, (\ref{haf}). Choose a \emph{compatible} almost
complex structure $J$ on $\C^n$, that is an almost complex
structure, $J^2 = -1$, so that $\omega_0(\cdot, J \cdot) =
g_{J}(\cdot, \cdot)$ is a Riemannian metric on $\C^n$. Then the
gradient lines of (\ref{haf}), satisfy Cauchy-Riemann type PDE,
\begin{equation}\label{CR}
\frac{\partial u}{\partial s} + J(t,u) \frac{\partial u}{\partial
t} + \nabla H(t,u) = 0 \textrm{  for  } u \in C^{\infty}( \R
\times S^1, \C ^n)
\end{equation}
Given two critical points $x_{+} , x_{-}$ of $\ah$, or
equivalently periodic orbits for the Hamilton equations $\dot{x} =
X_{H}(x); x(0)=x(1)$, consider the set of solutions,
$\mathcal{M}(x_{-} ,x_{+},J,H)$, of (\ref{CR}), such that,
$\lim_{s \to \pm \infty} u(s,t) = x_{\pm} (t) \}$. An element of
$\mathcal{M}(x_{-} ,x_{+},J,H)$ is called a \emph{Floer
trajectory}. The difference of the actions between the ends of a
Floer trajectory, $u$, is given by \emph{the energy}, $E_J(u)$,
\begin{equation}\label{eq14}
\ah ( x_{+}) - \ah (x_{-}) = \int_{\rs} g_J (\frac{\partial
u}{\partial s}, \frac{\partial u}{\partial s}) ds dt \equiv E_J
(u)\geq 0
\end{equation}
One can assign to each critical point, $x$ of $\ah$, under some
non-degeneracy assumptions, an index, $\mu _{CZ}(x)$, called
\emph{the Conley-Zehnder index}, cf. \cite{SZ}. Studying the
combinatorics of the solutions of (\ref{CR}), yields the
\emph{Floer homology groups}, $HF_{k} ^{[a,b)} (H,J)$, which are
independent of $J$ for generic choice of the almost complex
structure. These groups consists of linear combinations with
$\mathbb{Z}_2$ coefficients of critical points, $x$ of $\ah$, with
$\mu _{CZ}(x)=k$, $a \leq \ah(x) <b$. Now using the assumptions on
$\mathcal{H}_{ad}(U)$ and the functoriality of the Floer homology,
one defines, following \cite{FH1}, the Floer homology groups
$HF_{k} ^{[a,b)} (H,J)$ for $H \in \mathcal{H}_{ad}(U)$. Then the
symplectic homology groups $S_{k}^{[a,b)}(U)$ are defined as
direct limits of the Floer homology groups. In other words for an
\emph{exhausting (cofinal)} family, $(H_{\lambda}, J_{\lambda})$,
for $\mathcal{H}_{ad}(U)$, which is 1-parameter family such that
if $K \in \mathcal{H}_{ad}(U)$, there exists $\lambda'$ so that
$H_{\lambda} \geq K$, whenever $\lambda > \lambda'$. Then,
\begin{equation}\label{SH}
S_{*}^{[a,b)}(U) = \lim_{\lambda \to \infty} HF_{*} ^{[a,b)}
(H_{\lambda},J_{\lambda})
\end{equation}
Recall from \cite{FHW},
\begin{lm} \label{shlm} The symplectic homology groups of an open ball of radius
$R$, $B_R = B^{2n}(R) \subset \C^n$, satisfy
$$S^{[a,b)} _n (B^{2n} (R)) = \Z _2 \textmd{ for }   a
\leq  0 < b \leq \pi R^2, \textmd{ and } 0 \textmd{ otherwise}.$$
$$S^{[a,b)} _{n+1} (B^{2n}(R)) = \Z _2 \textmd{ for } 0 <a \leq \pi R^2 <
b, \textmd{ and } 0 \textmd{ otherwise}.$$
$$S^{[a,b)} _{k} (B^{2n}(R)) = 0 \textmd{ for } k<n \textmd{ or }
n<k <3n$$
\end{lm}
Let $U$ be an open and bounded subset of $\C^n$. Without loss of
generality we may assume that the origin, $0 \in U$. Let $r
> 0$ be a number such that $ B^{2n}(r) \subset U$. Pick numbers
$\varepsilon >0 $ such that $\varepsilon < \pi r^2$ and a number
$b > \pi r^2$. Observe that for large $b$, the natural map,
$$\Z _2 = S^{[0, \varepsilon)} _n ( B^{2n}(\rho)) \to S^{[0, b)}
_n (B^{2n} (\rho))$$ vanishes, (see \cite{V}). Let $R$ be
sufficiently large so that $B_r = B^{2n}(r) \subset U \subset
B^{2n}(R)= B_R$, then we have
$$\begin{CD}
\Z_2 = S^{[0, \varepsilon)} _n (B_R) @> i_R >> S^{[0,
\varepsilon)} _n(U)@> i_r >> S^{[0, \varepsilon)} _n (B_r) = \Z_2
\end{CD}$$
Where $i_r, i_R$ are the inclusion morphisms, cf. \cite{He, He1}.
Since the composition $i_R \circ i_r$ is an isomorphism, it
follows that $0 \neq \alpha_U = i_R(1) \in S^{[0, \varepsilon)} _n
(U)$. One then considers the natural map
$$i^b _U : S^{[0, \varepsilon)} _n (U) \to S^{[0, b)} _n (U)$$
and the Floer-Hofer (homological) capacity is defined as
\begin{equation}\label{fh}
c(U) = \inf \{ b \mid i^b _U (\alpha_U ) = 0 \}
\end{equation}

Next we present an alternative definition of the Floer-Hofer
capacity utilizing the action selector method. We follow
\cite{He1}, Section 3.2. Let $x_0 \in U$ and $B_r (x_0) \subset
U$, where $B_r(x_0)$ is the open ball centered at $x_0$ with
radius $r$. Denote by $\mathcal{H}_{ad}^{x_0, r}(U)$ the subset of
$\mathcal{H}_{ad}(U)$, consisting of functions $H$ with unique
minimum at $x_0$ so that $0> \min H=H(x_0)> -\pi r^2$. Now given
function $H \in \mathcal{H}_{ad}^{x_0, r}(U)$, we consider a
function $K \in \mathcal{H}_{ad}(B_r (x_0))$, so that $\min K =
K(x_0) = \delta < \pi r^2$ and $H \geq K$ on $S^1 \times \C^n$.
For $\delta <\varepsilon <\pi r^2$, we have that
$$HF_n ^{[0,\varepsilon)} (K) \simeq \mathbb{Z}_2$$
and the generator is $x_0$. Consider the monotonicity morphism,
\cite{FH1}, $$m : HF_n ^{[0,\varepsilon)} (K) \to HF_n
^{[0,\varepsilon)} (H)$$ and denote by $\alpha_H = m(x_0)$.
Consider the natural inclusion map, $$i^b : HF_n
^{[0,\varepsilon)} (H) \to HF_n ^{[0,b)} (H)$$ and define the
number
\begin{equation}\label{selec}
\sigma(H) = \inf \{b | i^b(\alpha_H)=0\}
\end{equation}
It is proven in \cite{He1}, Proposition 3.3 that $\sigma (H)$ is a
positive critical value of $\ah$.
\begin{rema} \label{remark} We remark that if $H \in \mathcal{H}_{ad}^{x_0,
r}(U)$ as above, it follows from the definition of Floer homology
groups that for some compatible almost complex structure $J$ there
exists a Floer trajectory for $(H,J)$, connecting $x_0$ and $x(t)$
where $x_t$ is a 1-periodic orbit of $X_H$, so that $\mu_{CZ}(x) =
n+1$ and $\sigma(H) =\ah (x)$.
\end{rema}
The action selector $\sigma$, defined above can be extended for
all $H \in \mathcal{H}_{ad}(U)$ and the Floer-Hofer capacity is
defined as
$$c(U) = \inf \{ \sigma(H) | H \in \mathcal{H}_{ad}(U)\}$$
We refer to \cite{He1}, for the details.
\section{Proof of Theorem \ref{D}.}
We recall the following lemma from \cite{B}, Lemma 1. Here we will
use the version for $\R^{2n}$ and $p=1$.
\begin{lm}\label{lm1} Let $U$ be a bounded domain in $\R^{2n}$ with
contact type boundary $S=\partial U$. There exists $\epsilon >0$,
an open neighborhood $V$ of $S$ in $\R^{2n}$ and a diffeomorphism
$\psi : S \times (-\epsilon,\epsilon) \to V$ so that:

i) For all $x \in S$ we have $\psi(x, 0) =x$;

ii) $\psi^* \omega = (1 + t) q^*(\omega_0|_{S}) + dt \wedge q^*(\alpha_0)$;\\
where $\alpha_0$ is a 1-form on $S$, such that, $d\alpha_0 =
\omega_0$ and $\alpha_0 \wedge \omega_0^{n-1} \neq 0$ on $S$, $q:
S \times (-\epsilon,\epsilon) \to S$ is the projection and $t$ is
a coordinate on $(-\epsilon,\epsilon)$.
\end{lm}
Set $$r=q\circ \psi^{-1} : V \to S$$
$$\beta_0 = r^{*} \alpha_0$$
$$\Omega = r^{*}(\omega_0|_{S})$$
$$z= t \circ \psi^{-1}$$
then
$$\omega_0 = (1+z) \Omega + dz \wedge \beta_0$$
Denote the Hamiltonin vector field corresponding to $z$, by $X_z$,
then we get
$$\beta_0 (X_z) =1$$
$$\Omega(X_z, \cdot) =0$$

Let $\epsilon >0$ be the number given by Lemma \ref{lm1}, we may
assume in addition that $1> \epsilon > 0$. Fix $\epsilon'$ such
that $\epsilon > \epsilon' >0$. For $0 <\tau \leq \epsilon$,
denote by $$V_{\tau} = \psi (N \times (-\tau, \tau)) = \{ x \in V|
z^2 (x) < \tau^2 \}$$ Consider the 1-form $B_0$ defined on $\C^n$
by $B_0 = f \beta_0$, where $f$ is a smooth function on $\C^n$
such that $f=1$ on $V_{\epsilon'}$ and $f = 0$ on $\R^{2n}
\setminus V_{\epsilon}$. This way we get one-form defined on
$\R^{2n}$ such that
\begin{equation} \label{B1}
B_0 = \beta_0 \textmd{ on } V_{\epsilon'}
\end{equation}
and
\begin{equation} \label{B2}
B_0 = 0 \textmd{ on } \C^n \setminus V_{\epsilon}
\end{equation}

Our next task is to build exhausting family for the set $U$.
Before we proceed we need, the following proposition.
\begin{prop}\label{b} The action spectrum $\Sigma (S)$ of a hypersurface
$S$ of contact type is a nowhere dense set.
\end{prop}
This proposition appears in many papers as a statement without a
proof. B. G\"{u}rel has carried out the details of the proof and I
thank her and V. Ginzburg for emailing me detailed outline of the
proof, \cite{Gu}.

Consider, for sufficiently large $\lambda \notin \Sigma(S)$,
smooth functions $g$ and $h$ on $\R ^{+}$ so that.
\begin{itemize}
\item $h'(t) = \lambda$ for $t \in [3 \lambda ^{-1}, 2 \lambda
^{-1} + \lambda^{-1/2}]$,

\item $h(t) = - \lambda ^{-1}$ for $t \in [0, \lambda ^{-1}]$,

\item $h(t) = - \lambda ^{-1} + \lambda^{1/2}$ for $t \geq 3
\lambda ^{-1} + \lambda ^{-1/2}$,

\item $h$ is convex on $[\lambda ^{-1}, 3\lambda^{-1}]$ and
concave on $[2 \lambda ^{-1} + \lambda ^{-1/2}, 3 \lambda ^{-1} +
\lambda ^{-1/2}]$,

\item $h (2 \lambda ^{-1}) <0$

\item $g(t) = - \lambda ^{-1} + \lambda^{1/2}$ for $ t <
(\lambda^{1/6} +1)^2 - \lambda^{-1}$,

\item $g'(t) = \mu/2$ for $t> (\lambda^{1/6} +1)^2$,

\item $g$ is convex on $ [(\lambda^{1/6} +1)^2 - \lambda^{-1},
(\lambda^{1/6} +1)^2]$.

\end{itemize}

Here $\mu \sim \lambda ^{1/6}$ and $\mu \notin \pi \Z$. We remark
that due to Proposition \ref{b}, we can perturb $h$ slightly, if
necessary, to assure that $h'(t) \notin \Sigma(S_t)$, for $t \in
[3 \lambda^{-1}, 2\lambda^{-1} + \lambda^{-1/2}]$. Now define
$H_{\lambda}$ as follows.
\begin{itemize}
\item $H_{\lambda} (x) = - \lambda ^{-1}$ for $x \in U$,

\item $H_{\lambda} (x) = h( z(x))$ for $x \in \bigcup_{0 \leq \nu
< 3 \lambda ^{-1} + \lambda ^{-1/2}} S_{\nu}$,

\item $H_{\lambda} (x) = g (|x|^2)$ for $ |x| > \lambda ^{1/6}$,

\item $H_{\lambda}(x) = - \lambda ^{-1} + \lambda^{1/2}$ for $x
\in B^{2n} (\lambda ^{1/6}) \setminus \{U \bigcup \bigcup_{0 \leq
\nu < 3\lambda ^{-1} + \lambda ^{-1/2}} S_{\nu}\}$.
\end{itemize}

It is clear that $\{H_{\lambda}\}$ is an exhausting family for
$U$. Next perturb each $H_{\lambda}(x)$, for $x \in U
\bigcup_{\tau \in [0,\lambda^{-1}]} S_{\tau}$, as in \cite{He,D},
to get non-degenerate family with unique prescribed global minimum
at the origin $0 \in U \setminus \bigcup_{\tau \in (-\epsilon,0)}
S_{\tau}$, which we assumed to be in $U$. Abusing the notation we
call the new family again $H_{\lambda}$. We can arrange the
perturbation in such a way, so that each $H_{\lambda}$, is bounded
together with its derivatives by $\lambda^{-1}$ on $\bar{U}$. In
particular we choose a constant $C_3 >0$ so that for $x \in
\bar{U}$ we have,
$$|X_{\hl}(x)| \leq C_3 \lambda^{-1}$$
 Next couple each $H_{\lambda}$ with a compatible non-degenerate almost complex
structure $J_{\lambda}$.

Standard arguments as in \cite{HZ1,He}, show that for sufficiently
large, but fixed $\lambda$, $0 < \sigma(H_{\lambda}) =
\mathcal{A}_{H_{\lambda}} (x_{\tau (\lambda)})$, where $x_{\tau
(\lambda)} \in S_{\tau (\lambda)}$ and $0 \leq \tau (\lambda) <
3\lambda ^{-1} + \lambda ^{-1/2}$. Moreover $x_{\tau (\lambda)}$ solves an equation
of the form, $\dot{x} = \rho X_z (x)$, where $\rho = h' (\tau
(\lambda))$. Our goal is to show that the length of $x_{\tau
(\lambda)}$ is bounded independently of $\lambda$.

We know from Remark \ref{remark}, that there is a Floer trajectory
for $\hl, J_{\lambda}$, $u$ s.t. $\lim_{s \to -\infty} u = 0$ and
$\lim_{s \to \infty} u = x_{\tau (\lambda)} (t)$. Then,
\begin{equation}\label{energy}
\mathcal{A}_{\hl} (x_{\tau(\lambda)}) - \mathcal{A}_{\hl}(0) =
E_{J_{\lambda}}(u) = \int_{\rs} \omega_0 (u_s, u_t -X_{\hl}(u))ds
dt
\end{equation}
Consider now holomorphic change of the variables $\sigma : \dot{D}
\to \R \times S^1$, where $\dot{D} = D\setminus \{0 \}$, and $D$
is the unit disk in $\C$. Set $v(s,t) =u(\sigma(s,t))$, here, by
abusing the notation we consider $s+ it$ to be the holomorphic
coordinate on $D$. Then we have that $v$, satisfies,
\begin{equation}\label{CR1}
\frac{\partial v}{\partial s} + J_{\lambda}(t,v) \frac{\partial
v}{\partial t} + \nabla H_{\lambda}(v) = 0
\end{equation}
with $\lim_{s \to 0} v(s,t) = 0$ and $\lim_{s \to 1} v(s,t) =
x_{\tau(\lambda)}(t)$. Moreover we have that, $$E_{J_{\lambda}}(u)
= E_{J_{\lambda}}(v)= \frac{1}{2}\int_{\dot{D}}
(|v_s|_{g_{J_{\lambda}}} ^2 + |v_t - X_{\hl}(v)|_{g_{J_{\lambda}}}
^2) ds dt$$ Write $\int_{\dot{D}} = \int_{\dot{D}} ^{I} +
\int_{\dot{D}}^{II} + \int_{\dot{D}}^{III} $, where
$\int_{\dot{D}} ^{I}$, means that the integration is taken over
$(s,t)$, for which $v(s,t) \in \bar{U}$, in $\int_{\dot{D}}^{II}$,
integration is over $(s,t)$ for which $v(s,t) \in \overline{U
\bigcup V_\epsilon} \setminus \bar{U}$, and in
$\int_{\dot{D}}^{III}$ the integration is over the remaining
values of $(s,t) \in \dot{D}$.

Consider the space of all almost complex structures $J$ on $\C^n$,
compatible with $\omega_0$. Denote, as before, by $g_J$ the
corresponding metric, i.e., $g_J(\cdot , \cdot) = \omega_0 (\cdot,
J \cdot)$. Since the set $\overline{U \bigcup V_\epsilon}$ is a
compact subset of $\C^n$, there is a constant $C_2 >0$ so that on
$\overline{U \bigcup V_\epsilon}$ we have that,
$$|\xi |_{g_J} \geq \sqrt{C_2} |\xi |_{g_{J_0}} = \sqrt{C_2}
|\xi|$$ for any $\xi \in \C^n$. Here $J_0 = i$ is the standard
complex structure on $\C^n$. Let $C_1 >0$ be a constant such that,
$$C_1 |\xi||\eta| \geq |dB_0(\xi,\eta)|$$
for all $ \xi, \eta \in \C^n$. Observe that we can extend
$X_{\hl}(v)$ smoothly over the puncture $s=0$ by setting it equal
to $0$. We have then,
\begin{align*}
\frac{1}{2}&\int_{\dot{D}}^{I} (|v_s|_{g_{J_{\lambda}}} ^2 + |v_t
- X_{\hl}(v)|_{g_{J_{\lambda}}} ^2) ds dt \geq \frac{C_2}{2}
\int_{\dot{D}}^{I} (|v_s| ^2 + |v_t - X_{\hl}(v)| ^2) ds dt\\
&\geq \frac{C_2}{2} \int_{\dot{D}}^{I} (|v_s| ^2  + |v_t|^2 /2-
|X_{\hl}(v)| ^2) ds dt \\&\geq C_2 /\sqrt{2} \int_{\dot{D}}^{I} |v_s||v_t|ds
dt - C_2/2 \int_{\dot{D}}^{I} |X_{\hl}(v)| ^2 dsdt \\& \geq
C_2/(\sqrt{2}C_1) \int_{\dot{D}}^{I} |dB_0(v_s,v_t)|ds dt - C_2/2
\int_{\dot{D}}^{I}C_3^2 \lambda^{-2} dsdt \\&\geq \frac{C_2}{\sqrt{2}C_1}
\int_{\dot{D}}^{I} |dB_0(v_s,v_t)|ds dt - \pi C_2 C_3^2
\lambda^{-2}/2
\end{align*}
On the other hand,
\begin{align*}
\frac{1}{2}&\int_{\dot{D}}^{II} (|v_s|_{g_{J_{\lambda}}} ^2 + |v_t
- X_{\hl}(v)|_{g_{J_{\lambda}}} ^2) ds dt \geq \frac{C_2}{2}
\int_{\dot{D}}^{I} (|v_s| ^2 + |v_t - X_{\hl}(v)| ^2) ds dt\\
&\geq \frac{C_2}{C_1} \int_{\dot{D}}^{II} |dB_0(v_s,v_t -
X_{\hl}(v) )| ds dt = \frac{C_2}{C_1} \int_{\dot{D}}^{II}
|dB_0(v_s,v_t)| ds dt
\end{align*}
The last inequality follows from the fact that
$dB_0(X_{\hl}(v),\cdot) =0$, when $v \in \overline{U \bigcup
V_\epsilon}$. Observe that $\int_{\dot{D}}^{III} |dB_0(v_s,v_t)|
ds dt = 0$. In view of the inequalities above we conclude that,
$$E_{J_{\lambda}}(v) \geq \frac{C_2}{\sqrt{2}C_1}
\int_{\dot{D}} |dB_0(v_s,v_t)|ds dt - \pi C_2 C_3^2
\lambda^{-2}/2$$ Further, applying Stokes' theorem, we get,
\begin{align*}
\int_{\dot{D}} |dB_0(v_s,v_t)|ds dt \geq |\int_{\dot{D}}
dB_0(v_s,v_t)ds dt| = |\int_{S^1} x_{\tau(\lambda)} ^* \beta_0| =
|h'(\tau(\lambda))|
\end{align*}
And we obtain that
\begin{align*}
\sigma(\hl) - \hl(0) = \mathcal{A}_{\hl}(x_{\tau(\lambda)}) -
\mathcal{A}_{\hl}(0)= E_{J_{\lambda}} (u) \geq
\frac{C_2}{\sqrt{2}C_1}|h'(\tau(\lambda))| - \pi C_2 C_3^2 \lambda^{-2}/2
\end{align*}
and this shows,
$$|\rho| = |h'(\tau(\lambda))| \leq \frac{\sqrt{2}C_1}{C_2}(c(\hl) - \hl(0) + \pi C_2 C_3^2
\lambda^{-2}/2)$$ Since $\sigma(\hl)$ is bounded independently of
$\lambda$, (it is less than the capacity of some large ball
containing $\overline{U \bigcup V_\epsilon}$), we conclude that
$|\rho|$ is bounded independently of $\lambda$, and therefore so is
the length of $x_{\tau(\lambda)}$. Applying Arzela - Ascoli theorem, we can find a
sequence $\{\lambda_k \to \infty\}$, s.t., $\lim_{k \to \infty}
\tau(\lambda_{k}) = 0$, $c(U) = \lim_{k \to \infty}
\sigma(H_{\lambda_k})$, and $x_0 = \lim_{k \to \infty}
x_{\tau(\lambda_{k})}$ is a closed characteristic on $S_0 = S$.
Moreover
$$c(U) = \lim_{k \to \infty}(\int_{S^1} x_{\tau(\lambda_k)} ^* \lambda_0 - \int_0 ^1 H_{\lambda_k}
(x_{\tau(\lambda_k)})dt) = \int_{S^1} x_0 ^* \lambda_0 \in
\Sigma(S)$$

%%%%%%%%%%%%%%%%%%%%%%%%%%%%%%%%%%%%%%%%%%%%%%%%%%%%%%%%%%%%%%%%%%%%%%%%%%%

%%%%%%%%%%%%%%%%%%%%%%%%%%%%%%%%%%%%%%%%%%%%%%%%%%%%%%%%%%%%%%%%%%%%%%%%%%%
%%%%%%%%%%%%%%%%%%%%%%%%%%%%%%%%%%%%%%%%%%%%%%%%%%%%%%%%%%%%%%%%%%%%%%%%%%%

\end{document}